\documentclass[11pt, english]{article}

\usepackage[margin= 1.7 cm]{geometry}
\usepackage[T1]{fontenc}
\usepackage{amsthm}
\usepackage{amsmath}
\usepackage{amssymb}
\usepackage{setspace}
\usepackage{mathtools}
\usepackage{verbatim}
\usepackage{booktabs}
\usepackage{graphicx}
\usepackage[hidelinks]{hyperref}
\usepackage{cleveref}
\usepackage{graphicx}
\usepackage{appendix}
\usepackage[inline]{enumitem}
\usepackage{framed}
\usepackage{subcaption}
\usepackage{microtype,xurl}
\usepackage{caption}
\usepackage{listings}
\usepackage{xcolor}

\usepackage{floatrow}
\usepackage[T1]{fontenc}

\lstdefinestyle{verifier}{
  language=Python,
  basicstyle=\ttfamily\fontsize{8.1}{10.2}\selectfont,
  keywordstyle=\bfseries,
  commentstyle=\itshape,
  stringstyle=\ttfamily,
  showstringspaces=false,
  keepspaces=true,
  columns=fullflexible,
  breaklines=true,
  breakatwhitespace=false,
  numbers=left,
  numberstyle=\tiny,
  numbersep=7pt,
  xleftmargin=1.6em,
  aboveskip=0.8em,
  belowskip=0.8em,
  tabsize=4,
  upquote=true
}

\theoremstyle{plain}

\newtheorem*{thm*}{Theorem}
\newtheorem{thm}{Theorem}
\Crefname{thm}{Theorem}{Theorems}

\newtheorem*{lem*}{Lemma}
\newtheorem{lem}[thm]{Lemma}
\Crefname{lem}{Lemma}{Lemmas}

\newtheorem*{claim*}{Claim}

\Crefname{claim}{Claim}{Claims}

\Crefname{prop}{Proposition}{Propositions}

\Crefname{cor}{Corollary}{Corollaries}

\Crefname{conj}{Conjecture}{Conjectures}

\Crefname{qn}{Question}{Questions}

\Crefname{obs}{Observation}{Observations}

\newtheorem{ex}[thm]{Example}
\Crefname{ex}{Example}{Examples}

\theoremstyle{definition}

\Crefname{prob}{Problem}{Problems}

\Crefname{defn}{Definition}{Definitions}

\newtheorem*{defn*}{Definition}

\theoremstyle{remark}

\usepackage{thmtools}
\usepackage{thm-restate}



\newcommand{\ceil}[1]{
    \left\lceil #1 \right\rceil
}

\newcommand{\eps}{\varepsilon}

\expandafter\def\expandafter\normalsize\expandafter{%
    \normalsize
    \setlength\abovedisplayskip{8pt}
    \setlength\belowdisplayskip{8pt}
    \setlength\abovedisplayshortskip{4pt}
    \setlength\belowdisplayshortskip{4pt}
}

\usepackage[square,sort,comma,numbers]{natbib}
\providecommand{\doi}[1]{\href{https://doi.org/#1}{doi:\nolinkurl{#1}}}

 \setlist[itemize]{leftmargin=*}

 \newcommand{\file}[1]{\mbox{\nolinkurl{#1}}}

\DeclareMathOperator{\tr}{tr}

\title{Polylogarithmic bound for the Gy\'arf\'as-Sumner conjecture}
\author{
Martin Bierbaumer\thanks{Faculty of Mathematics, University of Vienna, Vienna, Austria.}
\and
Matija Buci\'c\thanks{Faculty of Mathematics, University of Vienna, Vienna, Austria. This research was funded in whole or in part by the Austrian Science Fund (FWF) 10.55776/AST2416125. Email: \href{mailto:matija.bucic@univie.ac.at} { \ttfamily matija.bucic@univie.ac.at}. } 
}

\date{}

\begin{document}
\maketitle
\begin{abstract}
    The classical Gy\'arf\'as-Sumner conjecture asserts that for any fixed forest $F$, any graph $G$ not containing $F$ as an induced subgraph, with bounded clique number $\omega(G)$, has bounded chromatic number $\chi(G)$. Recently, Nguyen, Scott, and Seymour proved that every such graph has chromatic number at most $|G|^{o(1)}$. We improve this bound to a polylogarithmic one $(\log |G|)^{O(1)}$.
\end{abstract}

\section{Introduction}
One of the most basic questions in graph theory asks when the chromatic number of a graph is controlled by its clique number. It is immediate that it is a lower bound, but when can one also upper bound the chromatic number by a function of the clique number? 
Sadly, the answer is in general no, since one of the most classical results in probabilistic combinatorics, due to Erd\H{o}s \cite{erdos-girth-chromatic}, established the existence of graphs with high girth (so in particular no $K_3$) with arbitrarily large chromatic number. This construction is highly random, and leads to a very natural question of whether assuming the graph we are working with is structured changes the answer. This question lies at the core of the theory of chi-boundedness, which has attracted considerable attention over the years with many interesting results; see \cite{chi-bounded-survey} for a recent survey on the topic. One of the most natural structural assumptions that one might make is to require the graph to not contain a fixed small graph $H$ as an induced subgraph. It is immediate from the above-mentioned high-girth, high-chromatic number construction that this can not work if $H$ contains a cycle. This leads us to the conjecture, independently raised by Gy\'arf\'as and Sumner in the late 1970s and early 1980s, that forbidding a tree (or equivalently a forest) does suffice. 

This beautiful conjecture remains open, despite a lot of work on the topic. Let us highlight one of the first results on the topic due to Gy\'arf\'as, which settled the conjecture for paths \cite{gyarfas}, introducing a technique now called the Gy\'arf\'as path argument that has proved extremely useful over the years since. Gy\'arf\'as, Szemer\'edi and Tuza proved the radius-two case for triangle-free graphs \cite{gyarfas-szemeredi-tuza}. The latter was subsequently extended to arbitrary bounded clique number by Kierstead and Penrice \cite{kierstead-penrice}. Further special families of trees have been settled since then, including brooms, subdivisions of stars, and several related classes; we direct the reader to \cite{chi-bounded-survey} for more details and further references. 

For a general tree, the best previous result is due to a recent paper by Nguyen, Scott, and Seymour \cite{NSS} who proved that a graph $G$ not containing a tree $T$ as an induced subgraph, with $\omega(G) \le k$, satisfies $$\chi(G) \le |G|^{o(1)},$$ 
where the $o(1)$ term in the exponent decays only as a small\footnote{The exponent is $-1/((r-1)(k-1))$, where $r$ is the radius of the tree.} negative power of $\log |G|$. 
We improve this bound to a polylogarithmic one.

\begin{restatable}{thm}{main}\label{thm:main}
    Let $T$ be a fixed tree and $k$ be a positive integer. There exists $C_{T,k}>0$ such that every graph $G$ with $|G| \ge 2$, not containing $T$ as an induced subgraph, and $\omega(G) \le k$, satisfies
    $$\chi(G) \le C_{T,k}\:(\log |G|)^{14k}.$$
\end{restatable}

\subsection{Proof overview}

Our proof will be by induction on the clique number, so we assume that any $T$- free graph with a smaller clique number already has a small chromatic number, and our goal is to show this holds for the current graph $G$ as well. We will only use this inductive assumption on the neighborhoods of our vertices (which naturally have a smaller clique number) to argue that they can be colored using few colors. 
To bound the chromatic number of the whole graph, we will show that it must contain a large independent set, remove it, and iterate. So, in what follows, our goal is to find a large independent set in a $T$-free graph with small clique number, and with each vertex neighborhood being $ q$-colorable for some small $q$. Our proof is naturally split into several steps.

\textbf{Step 1.} We use $T$-freeness to find many closed walks of length $4$. This is done by observing that we can find many not necessarily induced copies of $T$, where many here is in terms of the density of $G$. Since each such copy must, in fact, not be an induced copy, it must contain an edge which completes a cycle with some subpath of $T$. This cycle is very short (since $T$ is of constant size), giving us a lot of short cycles in the graph, even when accounting for possible overcounting. Sidorenko-style inequalities allow us to translate this count to the desired count of length four closed walks. In reality, this step is not done via counting due to regularity issues, but is instead carried out via a stationary $T$-indexed random walk which attempts to embed a tree in a certain carefully chosen random way which allows us to do a ``weighted'' probabilistic version of the above counting argument; see \Cref{sec:T-free-to-walks} for more details.

\textbf{Step 2.} We use the high count of closed walks of length four to obtain many dense induced bipartite blocks with little density between them.
This is done by repeatedly extracting such blocks and finally subsampling them to ensure lower density between them. Here, the blocks themselves are provided by the high count of closed walks of length four by simply picking the edge which is used in the largest number of the walks and observing that each walk gives an edge between the neighborhoods of the endpoints of this edge. We can then invoke the fact that these neighborhoods are $ q$-colorable to find independent sets on both sides which still retain a lot of the density, giving us the desired block. We then remove this block and repeat so long as we have many closed walks. The final step, to ensure there is little density between the blocks, is to subsample the blocks with appropriate probability. We refer to this final collection of dense induced bipartite blocks with little density between them as a block-system. See \Cref{sec:walks-to-blocks} for more details.

\textbf{Step 3.} We use the block system to help us embed $T$. In fact, we need an additional assumption for now, which says that for any given block and any vertex in it, we can pick a large number of  candidates for its neighbors. These neighborhood candidate sets need to be reasonably large and well-behaved. The general strategy here is to fix a root of $T$, embed it arbitrarily, and then embed iteratively, by always picking a currently embedded vertex at even distance to the root which does not yet have all its children embedded, and then embed an additional child of that vertex in the same block and all its children (so grandchildren of the original one) in different previously unused blocks (by means of the candidate sets). So every even-distance vertex is embedded in a different block and is a unique vertex on its side of the bipartition, and all its children are embedded into the same block with it (on the other side).
For this process to work, we need several natural properties of the candidate sets; there is some flexibility here, but we go with: 
\begin{enumerate}[label=\alph*)]
    \item the neighborhood candidate sets of a vertex in block $i$ should not contain other vertices from block $i$, 
    \item they should intersect each block in at most one vertex, and \item neighborhood candidate sets of vertices in the same block should be disjoint.
\end{enumerate}

\textbf{Step 4.} The final part of the argument is to show that if we can never find suitable candidate sets (in a certain hereditary sense), then we get a large independent set. The existence of the candidate sets can conveniently be encoded by the existence of a large flow in a certain network, and by the max-flow min-cut theorem, this failure implies the existence of a small cut. 
This cut gives us an independent set with, in a certain sense, a small neighborhood. So we can iteratively remove these independent sets together with their neighbors to ensure that we can finally put all of them together into one large independent set. In reality, one needs to be a bit careful here, both with the fact that we can not remove all neighbors immediately (since what we actually know is that we have few neighbors outside of a few blocks, where we might have many neighbors) and we need to do some cleanup to maintain the properties of the block system which allow us to embed by Step 3, unless we again get the desired small cut. A final random selection of blocks and one bipartition side within each block then removes the remaining conflicts, yielding a large independent set. 

\section{Proof}

\subsection{Preliminaries}
In this section we introduce several preliminary lemmas and notation. 

\textbf{Notation.} Given a graph $G$ we denote by $V(G)$ and $E(G)$ its vertex and edge set respectively. We write $|G|=|V(G)|$. We denote by $e(G),\chi(G),\omega(G), \alpha(G)$ the number of edges, the chromatic number, the clique number, and the independence number of $G$ respectively. We denote by $\Delta_G,d_G,\delta_G$ the maximum, average and minimum degree of a graph $G$. We will sometimes also use $\Delta(G),d(G),\delta(G)$ for the same parameters to avoid double subscripts. Given $v \in V(G)$ we write $N_G(v)$ for the set of vertices adjacent to $v$ and $d_G(v)=|N_G(v)|$ for the degree of $v$ in $G$. 
We denote by $A_G$ the adjacency matrix of $G$. 
A graph \emph{homomorphism} from a graph $T$ into a graph $H$ is a function $\phi :V(T) \to V(H)$ such that whenever $\{v,u\} \in E(T)$ we also have $\{\phi(v),\phi(u)\} \in E(H)$. All our logarithms are in base $2$ unless otherwise specified.

The following is a standard induced regularization lemma; see \cite{regularization,regularization-old} for related statements.  It says that in an $n$-vertex graph we can get an induced subgraph with average and maximum degree only a logarithmic factor away, while still retaining at least half of the vertices. 
\begin{lem}\label{l:1}
    Any graph $G$ has an induced subgraph $H$ with $\Delta (H) \le 4 d(H) \log |G| $ with at least $|G|/2$ vertices.
\end{lem}

\begin{proof}
    Let $H_0=G$. For $i \ge 1$, given $H_{i-1}$, if $\Delta(H_{i-1})>4\log|G| \cdot d(H_{i-1})$ we let $H_i$ be the graph obtained from $H_{i-1}$ by removing all vertices with degree at least $\Delta(H_{i-1}) / 2$ in $H_{i-1}$. Otherwise, we stop the process and let $H=H_k$ be the final graph. Since in each step of this process the maximum degree gets halved, we know that the number of steps $k \le \log |G|.$ At any step $1\le i\le k$ we delete at most $\frac{|H_{i-1}|d(H_{i-1})}{\Delta(H_{i-1})/2}<\frac{|G|}{2\log|G|}$ vertices. Hence, $|H_k| \ge |G|/2$. 
\end{proof}

The next auxiliary lemma gives us an upper bound for the probability that the stationary random walk on a graph returns to a neighbor of its starting vertex, in terms of the number of closed walks of length four in the graph (which equals the trace of the fourth power of its adjacency matrix). Before stating it, let us introduce some background and notation. 

For a graph $G$ with no isolated vertices, we write $P$ for the transition matrix of the standard random walk on $G$, so that
\[
    P(u,v)=
    \begin{cases}
        1/d_G(u), & uv\in E(G),\\
        0, & \text{otherwise}.
    \end{cases}
\]
Its stationary distribution is $\pi(v)=d_G(v)/(2e(G))$. We will use that the walk is reversible, that is,
$\pi(u)P(u,v)=\pi(v)P(v,u)$ for all $u,v\in V(G)$, and hence a stationary random walk has the same distribution when read forwards or backwards. We refer to \cite{levin-peres-book} for further background on finite Markov chains and random walks on graphs.

\begin{lem}\label{lem:return-probability}
    Let $H$ be a graph with $\delta=\delta(H)>0$ and let $A$ be its adjacency matrix. Let $\ell \ge 2$, and let $X_0,X_1,\ldots, X_\ell$ be the stationary random walk on $H$. Then, $$\mathbb{P}(X_0 \sim X_\ell) \le \frac{\sqrt{|H| \tr A^4}}{2e(H)\sqrt{\delta}}.$$
\end{lem}
\begin{proof}
    Since we will work exclusively with $H$ in this proof, we omit it from all graph parameters. 
    
    Let $P$ be the transition matrix of our stationary random walk. 
    Observe that for any vertices $u,v$ and integer $i \ge 1$ we have 
    \begin{equation}\label{eq:uppr-bnd-transition}
        P^i(u,v) = \sum_{w} P^{i-1}(u,w)P(w,v)\le \frac1{\delta} \sum_{w} P^{i-1}(u,w)=\frac1{\delta}.
    \end{equation}
    
    In order to bound $\mathbb{P}(X_0 \sim X_\ell)$, it will be convenient to reveal $X_1$ and consider the probability that $X_0$ (which is a uniformly random neighbor of $X_1$ since the walk is stationary and reversible) and $X_\ell$ (which is simply the result of the random walk of length $\ell-1$ starting at $X_1$) are adjacent. To do so, let us fix a vertex $u$ in $H$, and we will condition our walk on $X_1=u$. Conditional on $X_1=u$ the variables $X_0$ and $X_\ell$ are independent, and for every vertex $v$ we have $\mathbb{P}(X_\ell=v \mid X_1=u)=P^{\ell-1}(u,v).$
    For any vertex $v$, the probability that a uniformly random neighbor of $u$ is also a neighbor of $v$ is equal to $|N(v) \cap N(u)|/d(u)=A^2(v,u)/d(u)$.
    So, 
    \begin{align*}
        \mathbb{P}(X_\ell \sim X_0 \mid X_1=u)&=\frac{1}{d(u)}\sum_{v \in V(H)} P^{\ell-1}(u,v) A^2(v,u)\\
        &\le \frac{1}{d(u)}\sqrt{\sum_{v \in V(H)} (P^{\ell-1}(u,v))^2}  \sqrt{\sum_{v \in V(H)} (A^2(v,u))^2}\\
        &\le  \frac{1}{d(u)} \cdot \frac{1}{\sqrt{\delta}} \cdot \sqrt{A^4(u,u)},
    \end{align*}
where we used Cauchy--Schwarz in the first inequality and \eqref{eq:uppr-bnd-transition} in the second to bound $\sum_{v \in V(H)} (P^{\ell-1}(u,v))^2\le \frac1{\delta}\sum_{v \in V(H)} P^{\ell-1}(u,v)=\frac1{\delta}$.

Finally, since the walk is stationary, we get 
$$\mathbb{P}(X_0 \sim X_\ell)=\sum_{u \in V(H)}\mathbb{P}(X_\ell \sim X_0 \mid X_1=u)\cdot \frac{d(u)}{2e(H)}\le \frac{1}{2e(H)\sqrt{\delta}} \sum_{u \in V(H)} \sqrt{A^4(u,u)}\le  \frac{\sqrt{|H| \tr A^4}}{2e(H)\sqrt{\delta}},$$
where the final inequality is another application of the Cauchy--Schwarz inequality.
\end{proof}

\textbf{Flows.} We use the standard terminology for flows in directed networks. A network consists of a directed graph with a source $s$, a sink $t$, and a nonnegative capacity on each edge. The value of a flow is the total flow leaving $s$, while the capacity of an $s$--$t$ cut is the sum of the capacities of the edges directed from the source side to the sink side. We use the max-flow min-cut theorem, which states that the maximum value of an $s$--$t$ flow equals the minimum capacity of an $s$--$t$ cut. We will also use the integrality theorem: if all capacities are integers, then there exists a maximum flow which is integral on every edge. See \cite{diestel} for further background.

\textbf{Ramsey.} We will also use the standard Erd\H{o}s-Szekeres \cite{erdos-szekeres} upper bound on the Ramsey numbers which guarantees that any graph with $\binom{k+t-2}{k-1} \le 2^{t+k-2}$ vertices contains either a clique on $k$ vertices or an independent set on $t$ vertices.

\subsection{From $T$-freeness to many short closed walks}\label{sec:T-free-to-walks}

The following is the key new lemma. It says that if a graph fails to have an induced copy of a fixed small tree $T$, then it needs to have a lot of short closed walks (specifically of length four). The intuition behind this is quite simple: one can attempt to run the greedy embedding argument starting from any vertex, and so long as we have some control on the minimum degree, we will get many, not necessarily induced, copies of $T$ in the graph. Since we do not have an induced copy of $T$, this implies that each of these embeddings must fail to be induced, giving rise to a short cycle supported on its vertex set. Accounting for how many times we count a specific cycle still results in a very large number of short cycles. Standard Sidorenko-style Cauchy-Schwarz inequalities then translate a high count of short cycles to a high count of ones of length four. In case the graph is regular, this becomes an easy counting argument. For irregular graphs, we would lose too much, and the following argument can be viewed as a certain weighted version\footnote{One can instead use a regularization lemma, and while this would in fact suffice to give a polylogarithmic bound it would introduce the size of the tree into the final exponent of the logarithm.}. The random embedding process we use below to deal with the irregularity is called a stationary $T$-indexed random walk; see \cite{odonnell-schramm,levin-peres} for more details on it, and \cite{conlon} for an application of it in a somewhat similar setting.

\begin{lem}\label{l:2}
    Let $T$ be a tree. Let $H$ be an $n$-vertex graph without $T$ as an induced subgraph. Then, the number of closed walks of length four in $H$ is at least $\frac{1}{16|T|^4} \cdot n d_H^3$.
\end{lem}

\begin{proof}
    If $d_H<4|T|^2$, then just by counting the length-four walks which return to the starting point both after two and four steps we get at least $\sum_{v\in H} d_H(v)^2\ge n d_H^2 \ge \frac{1}{4|T|^2}\cdot nd_H^3$, as desired. So let us assume $d_H \ge 4|T|^2$. Next, let us repeatedly delete any vertex with degree less than $d_H/4$ to obtain a graph $G$. Note that this deletes at most $n d_H/4=e(H)/2$ edges, so we obtain a graph $G$ with $N \le n$ vertices, with $\delta_G \ge d_H/4 \ge |T|^2$ and $e(G) \ge e(H)/2$. We proceed to count the number of closed walks of length four in $G$.
    
    We start by defining a random graph homomorphism $\phi:V(T) \to V(G)$ as follows. Let us fix an arbitrary vertex $u$ of $T$ as the root, and choose where to embed it according to the stationary distribution of the standard random walk on $G$. Namely, $\mathbb{P}(\phi(u)=v)=\frac{d_G(v)}{2e(G)}=:\pi(v)$, for every $v \in V(G)$. We proceed to recursively embed the remaining vertices of $T$ by embedding each vertex of $T$ uniformly at random to a neighbor of its parent. This gives us a probability distribution over homomorphic (but not necessarily injective) copies of $T$ in $G$. The probability of obtaining a fixed homomorphism $\eta: T \rightarrow G$ is 
    $$\pi(\eta(u))\left(\frac{1}{d_G(\eta(u))}\right)^{d_T(u)}\cdot \prod_{v \in T \setminus \{u\}}\left(\frac1{d_G(\eta(v))}\right)^{d_T(v)-1}=\frac{1}{2e(G)}\prod_{v \in T} d_G(\eta(v))^{1 - d_T(v)},$$
    where the first term comes from the probability of embedding the root at $\eta(u)$, the second from the probability of embedding each of its $d_T(u)$ children at the vertex specified by $\eta$, which happens with probability $1/d_G(\eta(u))$ per vertex, independently between vertices. Each term in the product comes from similar considerations for the remaining vertices $v \in T \setminus \{u\}$, except noting that they each only have $d_T(v)-1$ children. Note, in particular, that the final expression is independent of the choice of the root.

    As $G$ does not have an induced copy of $T$ as a subgraph, one of two things needs to happen: $\phi$ might not be injective, or two non-adjacent vertices of $T$ get mapped to adjacent vertices in $G$. In both cases, for some subpath of length at least $2$ from some vertex $a$ to some vertex $b$ in $T$, we embedded $b$ to the same vertex of $G$ as $a$ (in the first case) or to one of its neighbors (in the second case). Let us now analyze the probability that this happens for the unique path between any two specific distinct, non-adjacent vertices $a,b$ of $T$ (with the plan to take a union bound over all such pairs later). Let us denote the ``bad'' event that this happens by $B_{a,b}$.
    
    Due to the above established fact that the probability of picking any given homomorphism is independent of the choice of the root, we may assume that we picked $a$ as the root of $T$. Now, if $v_0=a,v_1,\ldots, v_\ell=b$ is the path from $a$ to $b$ in $T$, we have that $\phi(v_0),\phi(v_1),\ldots, \phi(v_\ell)$ traces a standard random walk on $G$, with the starting vertex being chosen from its stationary distribution. 
    So the probability that our bad event happens is precisely the probability that this random walk ends at a vertex of distance at most one from its starting point.

    The probability of return to the same vertex is at most $1/\delta_G$, since in the final step of the walk it had at least $\delta_G$ choices. The probability of ending in a neighbor of the starting vertex is upper bounded by \Cref{lem:return-probability}. So, if we denote by $W=\tr A_G^4$ the number of closed walks of length four in $G$ we get an upper bound 
    $\mathbb{P}(B_{a,b})\le \frac{1}{\delta_G}+\frac{\sqrt{|G| W}}{2e(G)\sqrt{\delta_G}}.$ There are at most $\binom{|T|}{2}$ choices for $a,b$, giving via a union bound
    $$1=\mathbb{P}\left(\bigcup_{a,b} B_{a,b}\right) \le \frac{|T|^2}{2\delta_G}+\frac{|T|^2\sqrt{|G| W}}{4e(G)\sqrt{\delta_G}},$$
    where the union is over all non-adjacent pairs of vertices of $T$ and the first equality holds since we know that one bad event \emph{must} happen.

    Since $\delta_G \ge |T|^2$ this implies 
    $$W \ge \frac{4(e(G))^2\delta_G}{|G||T|^4} \ge  \frac{(e(H))^2d_H}{4|H||T|^4}=\frac{1}{16|T|^4} \cdot |H| d_H^3,$$
    as desired.
\end{proof}

By taking $H$ to be a disjoint union of $K_{d+1}$'s, we get a graph with $\Theta(nd^3)$ closed walks of length four, which is clearly $T$-free for any $T$ such that $|T| \ge 3$, showing that the above lemma is tight in terms of its dependence on $|H|$ and $d(H)$.

\subsection{From many short walks to many useful blocks}\label{sec:walks-to-blocks}
The following lemma allows us to translate a good lower bound on the number of closed walks of length four to a special structure in the graph which in turn is going to be very helpful during the embedding process. The structure we aim to find is a collection of vertex-disjoint, dense, bipartite induced subgraphs, which we will refer to as blocks, with very few edges between them. This structure is inspired by the notion of a blockade, which is used heavily in the recent remarkable progress on the Erd\H{o}s-Hajnal conjecture; see  \cite{induced-density-i,induced-density-vi,induced-density-vii,erdos-hajnal-c5,pure-pairs-vi,pure-pairs-x} for just some examples and more details. 

\begin{defn*}
    A $D$-$(\lambda,\eps)$-\emph{block system} in a graph $G$ is a collection of disjoint $B_1,\ldots, B_I \subseteq V(G)$ such that 
    \begin{enumerate}
        \item\label{itm:1} $G[B_1 \cup \ldots \cup B_I]$ has maximum degree at most $D$.
        \item\label{itm:2} Each $B_i$ consists of at most $D$ vertices, and $G[B_i]$ is bipartite with minimum degree at least $\lambda D$.
        \item\label{itm:3} Every vertex of $B_i$ has at most $\eps D$ neighbors in any other $B_j$.
    \end{enumerate}
    We refer to $|B_1 \cup \ldots \cup B_I|$ as the \emph{volume} of the block system.
\end{defn*}

The following lemma allows us to extract a block system from a graph with many closed walks of length four. So let us set up some convenient notation first.
Given a graph $G$ with at least one edge we denote by $$\tau_4(G):=\frac{\tr A_G^4}{|G|d(G)\Delta(G)^2}.$$ 
Thus, $\tau_4(G)$ is the normalized count of closed walks of length $4$. For each vertex $v$, there are at most $d_G(v) \Delta(G)^2$ such walks starting at $v$ so $$\tr A_G^4 \le \sum_{v \in V(G)} d_G(v) \Delta(G)^2=|G| d(G) \Delta(G)^2,$$ so $\tau_4(G) \le 1$. On the other hand, \Cref{l:2} tells us that a $T$-free graph has $\tau_4(G) \ge \frac{1}{16|T|^4} \cdot \left(\frac{d(G)}{\Delta(G)}\right)^2.$

\begin{lem}\label{l:3}
    Let $G$ be a graph in which every neighborhood of a vertex is $q$-colorable, with $1 \le q \le \Delta(G)$. 
    For any $0<\eps\le 1$ we can find a  $D$-$(\lambda,\eps)$-block system $B_1,\ldots, B_I$ in $G$ with $D:=2\ceil{\frac{\Delta(G)}q}$, $\lambda=\frac{\tau_4(G)}{2^{11}}$, and 
    $$|B_1 \cup \ldots \cup B_I| \ge \frac{\eps\tau_4(G)^2d(G)}{2^{18}q\Delta(G)}\cdot |G|.$$
\end{lem}

\begin{proof}Throughout the proof let $n=|G|, \Delta=\Delta(G), d=d(G), \delta=\delta(G), \tau=\tau_4(G)$, so $\tr A_G^4=\tau n d \Delta^2$.

We first show that so long as we haven't deleted too many vertices, we can always extract another bipartite dense block (with property \ref{itm:2}). To this end, note that deleting a vertex deletes at most $4\Delta^3$ closed walks of length four. So, as long as fewer than $\frac{\tau n}{8} \cdot \frac{d}{\Delta}$ vertices have been removed, the remaining graph $G'$ still has $\tr A_{G'}^4 \ge \frac12 \tau n d \Delta^2$ closed walks of length four. Double counting these walks and edges contained in them, we can find an edge $uv$ of $G'$ such that there are at least $\tau \Delta^2/2$ of these walks of the form $u,a,b,v,u$. This gives us at least $\tau \Delta^2/2$ ordered pairs of adjacent vertices $(a,b)$, with $a \in N_{G'}(u)$ and $b \in N_{G'}(v)$.
This gives us disjoint sets $U \subseteq N_{G'}(u)$ and $V \subseteq N_{G'}(v)$ with at least $\tau \Delta^2/8$ edges between them, by the standard first moment argument where we assign any vertex in the intersection to only one of $U$ or $V$ uniformly at random.

Since both $U$ and $V$ are subsets of vertex neighborhoods, we can partition them into at most $q$ independent sets. Now further partition these independent sets into as few pieces of size at most $\ceil{\Delta/q}=D/2$ as possible. This results in up to $2q$ independent subsets on each side, so we can find independent sets $U' \subseteq U$ and $V' \subseteq V$ each of size at most $D/2$ with at least $\frac{\tau \Delta^2}{32q^2}\ge \tau D^2/2^{9}$ edges between them. Finally, repeatedly delete any vertex with degree less than $\tau D/2^{10}$, which can delete at most half of the edges. This results in an induced bipartite subgraph with at most $D$ vertices, minimum degree at least $\tau D/2^{10}$ and maximum degree at most $D/2$. So in particular, it satisfies our desired conditions in property \ref{itm:2}.

Let $A_1,\ldots, A_J$ be a maximal vertex-disjoint collection of such bipartite induced subgraphs (which we will refer to as blocks). By the discussion above, this family contains at least $\frac{\tau n d}{8\Delta}$ vertices. We now subsample these bipartite subgraphs to ensure the sparsity condition (property \ref{itm:3}) between different blocks. 
With this in mind, let us activate the blocks independently, each with probability $p:=\frac{\eps \tau }{2^{14}q}\le 1$.

For any vertex $v \in A_i$, there are at most $\Delta/(\eps D) \le q/\eps$ blocks $A_j$, with $j \neq i$, in which $v$ has more than $\eps D$ neighbors. We say $v$ is bad if one of these blocks activates. So by a union bound $v$ is bad with probability at most $pq/\eps$. We additionally declare $v$ bad if it has more than $D/2$ neighbors across other activated blocks. Since the expected number of neighbors of $v$ in activated blocks is at most $p \Delta$, Markov's inequality tells us that this happens with probability at most $2p\Delta/D \le pq$, so in total $v$ is bad with probability at most $pq/\eps+pq \le 2pq/\eps \le 1/4.$ 

Let us now discard any activated block with more than $\tau D/2^{11}$ bad vertices entirely, and remove any bad vertices from any other block. Since the expected number of bad vertices in a block is at most $D \cdot 2pq/\eps$,  another application of Markov's inequality implies that the probability that a block is discarded is at most $\frac{2^{12}pq}{\eps \tau}= \frac14$. The badness of vertices in $A_i$ depends only on the activation of other blocks, so both estimates remain valid conditional on $A_i$ being activated. In total, the probability that a vertex $v$ belongs to an activated block which was not discarded and $v$ is not bad is at least $p/2$. Let us now fix an outcome in which at least $\frac{p}2\cdot \frac{\tau n d}{8\Delta} = \frac{\eps \tau^2 n d}{2^{18}q\Delta}$ vertices survive and call the remaining parts of nonempty blocks $B_1,\ldots,B_I$. 

Since we deleted all bad vertices, we ensure that $G[B_1 \cup \ldots \cup B_I]$ has maximum degree at most $D/2+D/2=D$ and that property \ref{itm:3} holds. Since from any block, we deleted at most $\tau D/2^{11}$ vertices, and we had minimum degree of $\tau D/2^{10}$ in each block, the minimum degree remains at least $\tau D/2^{11}$ in each $G[B_i]$, so property \ref{itm:2} still holds.
\end{proof}

\subsection{Using the block system to embed}
The following definition introduces for every vertex $v$ of a block system a set $Y_v$ which one should think of as the candidate set for where neighbors of $v$ can be embedded. We impose several restrictions on these candidate sets in order to help with the embedding.
\begin{defn*}
    Given a block system $B_1,\ldots, B_I$ in a graph $G$, write $B=B_1 \cup \ldots \cup B_I$. We call a collection of sets $(Y_v)_{v \in B_1 \cup \ldots \cup B_I}$ a \emph{neighborhood candidate system} if the following hold. 
    \begin{enumerate}[label=\alph*)]
        \item\label{itm:ncs-1} If $v\in B_i$ then $Y_v$ is a subset of the neighbors of $v$ in $B \setminus B_i$.
        \item\label{itm:ncs-2} Each $Y_v$ intersects any block in at most one vertex.
        \item\label{itm:ncs-3} $Y_v$ and $Y_u$ are disjoint for any distinct $v,u$ in the same block $B_i$.
    \end{enumerate}
\end{defn*}

The following lemma says that if we can find such a candidate system, where each of the candidate sets is large then we can embed an induced copy of any tree $T$. We note that \Cref{l:3} will provide us with a block system, but not necessarily a suitably large neighborhood candidate system. Instead, we will see in the final lemma in the following subsection that if we fail to find it then we can find a large independent set. 

\begin{lem}\label{l:5}
    Fix a tree $T$ and $k \in \mathbb{N}$. 
    Suppose that for some $D>0$ and $4|T|/D\le \lambda \le 1$ we are given a non-empty $D$-$(\lambda,\lambda/(4|T|))$-block system $B_1,\ldots, B_I$ in a graph $G$ with $\omega(G) \le k$ and a neighborhood candidate system $(Y_v)_{v \in B_1 \cup \ldots \cup B_I}$ such that $|Y_v|\ge {2^{|T|+k+2}}/{\lambda}$ for each $v$. Then $G$ contains $T$ as an induced subgraph.
\end{lem}

\begin{proof}
    Write $B=B_1 \cup \ldots \cup B_I$. Since all embedding choices lie in $B$, we may replace $G$ by $G[B]$, and hence assume that $\Delta(G) \le D$. 

    Choose an arbitrary vertex $r$ as the root of $T$. We will embed $T$ starting from the root and embedding vertices layer by layer in order increasing in distance from the root. In particular, we will always embed vertices in an even layer into a new block (using the neighborhood candidate system), while children of all vertices in even layers will be embedded inside the same block as their parent (using the high density assumption inside the blocks and the bipartite structure). Thus, each used block contains at most one vertex from an even layer together with some of its children.

    Embed the root $r$ at an arbitrary vertex of $B_1$. 
    At each step, given an already embedded vertex $p$ of $T$ at even distance from $r$, for which we have not yet embedded all of its children, we will embed one such new child $c$ of $p$ into the same block $B_i$ where $p$ was embedded, and all $c$'s children into new previously unused blocks. Let $a \in B_i$ be the vertex where we embedded $p$. Vertex $a$ has at least $\lambda D$ neighbors in $B_i$. Each of the other, at most $|T|$, already embedded vertices of $T$ has at most $\lambda D/(4|T|)$ neighbors in the bipartition side of $B_i$ not containing $a$ (the other children of $p$ are the only vertices embedded inside $B_i$ and they are all embedded on the same side of the bipartition, the one not containing $a$). In addition, we have already embedded at most $|T|$ other children of $p$ into $B_i$ so in total there are at least $\lambda D-\lambda D/4-|T|\ge \lambda D/2$ candidates for $c$, so unused and with no edges to any previously embedded vertex (except $p$).

    Let $F$ be the union of already used blocks and neighborhoods of already embedded vertices. So, in particular $|F| \le 2|T|D$. Since the sets $Y_v$ with $v$ ranging across the at least $\lambda D/2$ candidates are disjoint (since the candidates all belong to the same block $B_i$) we know that there exists a candidate vertex $v$ with $|Y_v \cap F|\le \frac{ 2D|T|}{\lambda D/2}=\frac{4|T|}{\lambda}$. In particular, $|Y_v \setminus F| \ge 2^{|T|+k}$, so by Ramsey's theorem, since $\omega(G) \le k$, we can find an independent set of size $|T|$ inside $Y_v \setminus F$. In particular, we can embed $c$ into $v$ and its children into this independent set, since $Y_v$ contains at most one vertex per block by property \ref{itm:ncs-2} we get that these children are embedded into different blocks. So we can continue the process until we embed the whole $T$.
\end{proof}

\subsection{Failure to embed results in a large independent set}

The following lemma shows that (repeated) failure to embed $T$, meaning non-existence of a suitable neighborhood candidate system by \Cref{l:5}, results in a large independent set. 

\begin{lem}\label{lem:failure-to-embed}
    Fix a tree $T$ and $k \in \mathbb{N}$. 
    Suppose that for some $D>0$ and $0< \lambda \le 1$ we are given a $D$-$(\lambda,\lambda/(8|T|))$-block system $B_1,\ldots, B_I$ in a $T$-free graph $G$ with $\omega(G) \le k$. Then, 
    $$\alpha(G) \ge \frac{\lambda^3}{2^{2|T|+2k+12}} \cdot |B_1 \cup \ldots \cup B_I|.$$ 
\end{lem}

\begin{proof}
Let $B=B_1 \cup \ldots \cup B_I$. If $\lambda D < 8|T|$, then by the standard greedy removal bound we have 
    $$\alpha(G) \ge \alpha(G[B]) \ge \frac{|B|}{\Delta(G[B])+1}\ge \frac{|B|}{D+1} \ge \frac{\lambda}{16|T|} \cdot |B|,$$
which is larger than the desired bound. So let us assume $\lambda D \ge 8|T|$.

Throughout the process, we will remove some vertices. At every stage, we will ensure that the remaining vertices are split into blocks which are subsets of the original $B_i$'s, with any block whose minimum degree falls under $\lambda D/2$ being completely removed. This ensures that at any point, we are working with a $D$-$(\lambda/2, \lambda/(8|T|))$-block system $B_1' \cup \ldots \cup B_{J}'$, possibly with $J < I$. Let $B'=B_1' \cup \ldots \cup B_J'$.

As long as $B'$ is not empty, the inequalities $1\ge \lambda/2 \ge 4|T|/D$, and the fact that $G$ is $T$-free, allow us to apply \Cref{l:5} to the current block system. Thus, no neighborhood candidate system with each candidate set of size at least $2^{|T|+k+3}/{\lambda}$ exists. Since the candidate system requirements for different source blocks do not constrain one another, there must be a current block $B_i'$ for which we cannot find sets $Y_v$, one for each $v \in B_i'$, each of size at least $2^{|T|+k+3}/\lambda$ and satisfying the three defining properties of the neighborhood candidate system. We encode this as the nonexistence of a certain integral flow. 

To this end, we construct an auxiliary (directed) flow network by adding a source vertex, three intermediate layers, and a sink vertex, where 
\begin{itemize}
    \item the first layer consists of vertices in $B_i'$,
    \item the second layer consists of pairs $(v,B_j')$ with $v\in B_i'$ and $j \neq i$, and 
    \item the third layer consists of $B' \setminus B_i'$.  
\end{itemize} 
We connect the source to every vertex in the first layer by a directed edge of capacity $\eta:=\ceil{2^{|T|+k+3}/\lambda}$. Every first layer vertex $v \in B_i'$ is connected to every vertex $(v,B_j')$ in the second layer by a directed edge of capacity $1$. Every vertex $(v,B_j')$ from the second layer is connected to every neighbor $w$ of $v$ in $B_j'$ in the third layer by a directed edge of infinite capacity (in reality we take a large finite capacity), and every vertex of the third layer is connected to the sink by a directed edge of unit capacity. The following illustration of the network and the capacities might help with visualization.
$$\text{source} \stackrel{\eta}{\longrightarrow} \underbrace{v}_{\in B_i'} \stackrel{1}{\longrightarrow} \underbrace{(v, B_j')}_{v \in B_i',\ j \neq i} \ \stackrel{\infty}{\longrightarrow} \underbrace{w}_{\in N(v) \cap B_j'} \stackrel{1}{\longrightarrow} \text{sink}.$$

Observe that if we can find an integral flow of value $\eta |B_i'|$ in this network, we can define $Y_v$ to be the set of $w$ in the third layer that receive a positive flow from a vertex $(v,B_j')$ for some $j$. Indeed, since there is a total capacity of $\eta |B_i'|$ from the source to the first layer, each $v \in B_i'$ sends exactly $\eta$ units of flow to the second layer. Since the flow is integral, this means there are exactly $\eta$ sets $B_j'$ with $j\neq i$ for which $v$ sends a unit flow to $(v,B_j')$. Now each vertex $w$ in the third layer $B' \setminus B_i'$ can receive at most one unit of flow (due to the capacity of its edge to the sink). This implies that exactly $\eta$ neighbors of $v$ in $B' \setminus B_i'$ get placed into $Y_v$, establishing property \ref{itm:ncs-1} and the requirement on the size of $Y_v$. The same condition ensures that each $w$ is assigned to at most one $Y_v$, establishing disjointness of these sets (property \ref{itm:ncs-3}).
The capacity of one from $v$ to the second layer ensures that only one vertex is chosen from each $B_j'$ with $j \neq i$, establishing property \ref{itm:ncs-2}. 
Since such a neighborhood candidate system does not exist, we conclude that such an integral flow cannot exist. 

Since in a network with all edge capacities integral there is always an integral maximum flow, we conclude by the max-flow min-cut theorem that there exists a cut in this network with capacity less than $\eta|B_i'|$. Let $X \subseteq B_i'$ and $Z \subseteq B' \setminus B_i'$ be the sets of vertices of the first layer and third layer, respectively, which belong to the source side of the cut. 
Since we can not have infinite-capacity edges crossing the cut, we know that for any $x \in X$ and $j \neq i$ for which $N_G(x) \cap B_j'$ contains a vertex outside of $Z$, we must have $(x, B_j')$ on the other side of the cut. So if we denote by $b_Z(x)$ the number of such blocks, the cut capacity is at least
$$\eta(|B_i'|-|X|)+ \sum_{x \in X} b_Z(x)+|Z|.$$
Since we know it is less than $\eta|B_i'|$ we conclude that 
\begin{equation}\label{eq:ind-set-expansion}
    |X|\eta > |Z|+\sum_{x \in X} b_Z(x).
\end{equation}
This inequality tells us that, in a certain sense, $X$ is large compared to its neighborhood outside of $B_i'$. More specifically, if we delete $Z$ from our graph (which we know isn't much larger than $X$ by the above inequality), the remaining neighbors of $X$ are on average concentrated in only a few blocks. 
So our strategy will be to place $X$ into a ``would be'' independent set $U$ and delete $X \cup Z$ from our block-system, remove any blocks whose minimum degree dropped below $\lambda D/2$, and then repeat until there are no vertices left in the block-system. 
At the end of the process, we obtain a set $U$, which, while not quite independent, is rather sparse in a certain well-specified way and will allow for an extraction of an actual rather large independent set. 

Let us first show that $U$ indeed captures a large portion of our original block-system vertex set $B$. Let $X_1,\ldots, X_s$ and $Z_1,\ldots, Z_s$ denote the sets $X$ and $Z$ which we deleted during the above process. So in particular, $U=X_1 \cup \ldots \cup X_s$ and $B$ is the union of $X_i$'s, $Z_i$'s and the blocks which were deleted entirely due to their minimum degree dropping below $\lambda D/2$. By \eqref{eq:ind-set-expansion}, we know that $|X_i| \ge \frac{1}{1+\eta} \cdot (|X_i|+|Z_i|)$. In order to control the blocks that were deleted due to the appearance of a vertex of degree less than $\lambda D/2$, notice that in order for such a vertex to appear, we needed to delete at least $\lambda D/2$ vertices from that block by deleting $X_i \cup Z_i$'s. In particular, since each block has size at most $D$ we know that $X_1 \cup \ldots \cup X_s \cup Z_1 \cup \ldots \cup Z_s$ captures at least a $\lambda/2$ proportion of the original vertices of the     blocks that were deleted. Putting all of these observations together, we conclude that 
\begin{equation}\label{eq:U-lwr-bnd}
    |U| \ge \frac{|B|}{(1+\eta)(1+2/\lambda)} \ge \frac{|B|\lambda}{4\eta}.
\end{equation}

We are now ready for the final cleanup step in order to convert $U$ into an actual independent set. To do so, we subsample each original block with probability $p=1/(2 \eta)$ and for each sampled block retain only one of its two sides uniformly at random. Our final independent set will consist only of vertices $x$ of $U$, say $x\in X_i$, which were retained after this process and for which none of the $b_{Z_i}(x)$ blocks containing neighbors in the $i$-the deletion step were retained. To see that this process indeed produces an independent set, notice that only keeping one side of the bipartition ensures that all $x$ within the same block are independent, while for retained vertices $x \in X_i$ and $x' \in X_j$ for $i<j$ which belong to different blocks note that if $x \sim x'$ then since $x'$ was not deleted as part of $Z_i$ we know that the block containing $x'$ would contain a neighbor of $x$ not in $Z_i$ and would by the rules of the above process not be retained. 

The final step is to compute the expected number of retained vertices from $U$. For every $x \in U$ let $b_x:=b_{Z_i}(x)$, where $x \in X_i$. So, any $x \in  U$ is retained with probability $p \cdot \frac12 \cdot (1-p)^{b_x}\ge \frac p2 \cdot (1-pb_x).$ Therefore, the expected size of our independent set is at least
$$\sum_{x \in U} \frac p2 \cdot (1-pb_x)=\frac p2 \cdot \left(|U|-p\cdot \sum_{x \in U} b_x\right)\ge \frac p2 \cdot (1-\eta p) |U| \ge \frac{|U|}{8 \eta},$$
where in the first inequality we used \eqref{eq:ind-set-expansion} to get $\sum_{x \in U} b_x = \sum_{i=1}^s \sum_{x \in X_i} b_x \le \sum_{i=1}^s \eta |X_i|=\eta |U|.$

Combining this with the above lower bound on $|U|$ we get
$$\alpha(G) \ge \frac{|U|}{8 \eta} \ge \frac{|B|\lambda}{32\eta^2}\ge \frac{|B|\lambda^3}{2^{2|T|+2k+12}},$$
as desired.
\end{proof}

\subsection{The final accounting}
We are now ready to put together a proof of \Cref{thm:main}.

\main* 

\begin{proof}
    The proof is by induction on $k$, with the case $k=1$ being immediate with $C_{T,1}=1$. Suppose now that the result holds for any $T$-free graph with clique number at most $k-1$.

    Let $G$ be an $n$-vertex $T$-free graph with $\omega(G) \le k$. By Lemma~\ref{l:1} we can find an induced subgraph $H$ with 
    $$|H| \ge \frac{n}2, \quad\quad {d(H)} \ge \frac{\Delta(H)}{4\log n}.$$
    If $H$ has no edges, it is itself an independent set larger than what we will establish below, so let us assume that $H$ has at least one edge.
    
    By Lemma~\ref{l:2} we have 
    $$\tau_4(H) \ge \frac{1}{16|T|^4} \cdot \left(\frac{d(H)}{\Delta(H)}\right)^2 \ge \frac{1}{2^{8}|T|^4 \log^2 n}.$$ 
    
    Since every neighborhood of a vertex induces a graph with clique number at most $k-1$, by induction we can color it using at most $q:=\min \{C_{T,k-1}(\log n)^{14(k-1)},\Delta(H)\}$ colors. Using Lemma~\ref{l:3} with this choice of $q$, $D=2\ceil{\Delta(H)/q}$, $\lambda= \tau_4(H)/2^{11}$, and $\eps = \lambda/(8|T|)$ we get a $D$-$(\lambda, \eps)$-block system in $H$ with volume at least 
    $$\frac{\eps \tau_4(H)^2 d(H)}{2^{18}q \Delta(H)} \cdot |H| \ge \frac{n}{C_{T,k-1}|T|^{13}2^{59}\log^{14k-7} n}.$$

    Combined with \Cref{lem:failure-to-embed} we find an independent set in $G$ of size at least 
    $$\frac{\lambda^3}{2^{2|T|+2k+12}} \cdot \frac{n}{C_{T,k-1}|T|^{13}2^{59}\log^{14k-7} n} \ge \frac{n}{C_{T,k}\log^{14k-1} n},$$
    where we set $C_{T,k}=C_{T,k-1}2^{2|T|+2k}|T|^{25}2^{128}$.
    
    Since both $T$-freeness and $\omega(G) \le k$ are hereditary properties we can remove this independent set and repeat. After $C_{T,k} \log^{14k} n$ steps we remove the whole graph showing that $\chi(G) \le C_{T,k} \log^{14k} n$, completing the induction and the proof. 
\end{proof}

By solving the recursion we obtained on $C_{T,k}$ we get that we can take $C_{T,k} \le 2^{(2+o(1))|T|k+k^2}.$

\section{Concluding remarks}\label{sec:conclusion}

The most natural remaining question is, of course, to remove even the remaining polylog factor from \Cref{thm:main} and prove the Gy\'arf\'as--Sumner conjecture in full. An intermediate, likely much easier goal might be to at least remove the dependence of the exponent of the polylog factor on $k$. 
We note that we made no particular effort in optimizing the constants in our result. There are several ways apparent to the authors to do this. For example, instead of using the stationary $T$-indexed random walk to get many different $C_4$'s in \Cref{l:2}, we can stop the argument short and simply extract high typical density between neighborhoods of vertices, which can be used directly to find a block system. This avoids going back and forth to the count of four cycles, which avoids several applications of Cauchy-Schwarz and results in a much better final exponent. We opted for including the above argument as we find it significantly more intuitive and since we find \Cref{l:2} interesting in its own right. The first author will include some of these optimizations in his master's thesis.   

In \cite{subdivisions}, Nguyen, Scott, and Seymour proved a bound of exactly the same polylogarithmic form on the chromatic number for graphs which do not contain an induced subdivision of $K_t$, for a fixed $t$. We note that there are some parallels between the arguments, namely their star systems play a similar role to our block systems, and we both follow the same overall strategy of removing independent sets which do not ``expand''. The notions of expansion we are working with are very different, and we in fact do allow our removed independent sets to expand, but in a certain controlled way which allows the final subsampling cleanup to go through. We note that this iterative removal of non-expanding independent sets is a standard idea in the area dating back at least to a paper \cite{ErdosFaudreeRousseauSchelp1978} of Erd\H{o}s, Faudree, Rousseau, and Schelp from 1978. Despite some similarities in the overall approach, we do not see how to obtain a unified argument that would capture both results.

Another natural direction would be to at least try to prove a linear lower bound on the size of the independent set. This avoids one of the seemingly major difficulties in proving the full conjecture, namely that both our argument and the argument of \cite{NSS} iteratively remove independent sets, which introduces an additional logarithmic factor loss for the chromatic number bound. 
We note that, similarly to \cite{NSS}, our argument can, with minimal modification, be used to prove a lower bound on the independence number of the form $|G|/ O_{|T|,k}(\log^{14k} \Delta(G))$. 

Our bound on closed walks (\Cref{l:2}) has an easy consequence for the additive structure of connection sets of Cayley graphs. Let $(\Gamma,+)$ be a finite abelian group and $\varnothing \neq S \subseteq \Gamma \setminus \{\textrm{id}\}$ with $S=-S$. The additive energy of $S$ is defined as $E(S):=|\{(a,b,c,d) \in S^4: a+b=c+d\}|$.
By an immediate application of \Cref{l:2} if the Cayley graph $\operatorname{Cay}_{\Gamma}(S)$ contains no induced copy of a tree $T$ with $|T| \ge 3$ then $E(S) \ge |S|^3/(16|T|^4)$. Consequently the Balog--Szemer\'edi--Gowers theorem, in the quantitative form of Reiher and Schoen~\cite{reiher-schoen-bsg} gives a subset $S' \subseteq S$ of size at least $\frac{c}{|T|^2}\cdot |S|$ with bounded doubling, namely $|S'-S'| \le C|T|^{16}|S'|$, where $c,C>0$ are absolute constants. Since sets with small doubling are often viewed as very structured in additive combinatorics this gives a curious connection between graph theoretic notion of structure and additive one. We note that the question of forbidding induced trees in Cayley graphs dates back to a paper of Babai and S\'os \cite{babai-sos} from 1985. Also to illustrate that this condition can be quite natural from the additive point of view consider $\operatorname{Cay}_{\Gamma}(S)$ and $T$ being a star on $r$ elements $K_{1,r}$. Then forbidding an induced $K_{1,r}$ translates to forbidding existence of $s_1,\ldots, s_r \in S$ for which the differences $s_i-s_j \notin S$ for every $i \neq j$. The more classical variant of this property is forbidding the sums (which is equivalent over characteristic two) bringing us to the famous Erd\H{o}s--Moser sum free problem, where much stronger results than our easy observation are known, see \cite{tao-vu-sum-avoiding,sanders-erdos-moser}.

\medskip 

\subsection*{Statement on the use of AI}
ChatGPT 6 Astra Pro was used to proofread and simplify the original version of the argument, which is due to the authors and which supplied a similar polylogarithmic result except with the power of the logarithm also depending on $|T|$. In particular, it suggested the use of a stationary $T$-indexed random walk, which replaced the original counting-plus-regularization argument due to the authors.
It also supplied the flow-based formulation of the argument behind \Cref{lem:failure-to-embed}, which significantly simplified an original iterative argument based on repeated use of Hall's Theorem.  

\medskip 
\subsection*{Acknowledgements} The second author would like to thank Lior Gishboliner for useful early discussions on this problem during his visit to the University of Toronto in April 2025.
He would also like to thank Tung Nguyen, Alex Scott, and Paul Seymour for useful discussions about related problems, which heavily inspired \Cref{l:3}, and Tung Nguyen for useful comments on an early draft.

\providecommand{\MR}[1]{}
\providecommand{\MRhref}[2]{%
\href{http://www.ams.org/mathscinet-getitem?mr=#1}{#2}
}

\bibliographystyle{amsplain_initials_nobysame}
\bibliography{gyarfas-sumner-references}

\end{document}